\documentclass[12pt]{article}
\usepackage{amssymb,amsmath}
\usepackage{graphicx}
\usepackage{amsfonts}
\usepackage{float}
\usepackage{flafter}
\usepackage{color}

\labelwidth\leftmargini\advance\labelwidth-\labelsep

\def\R{\relax\ifmmode I\!\!R\else$I\!\!R$\fi}

\def\Z{\relax\ifmmode Z\!\!\!Z\else$Z\!\!\!Z$\fi}

\def\C{\relax\ifmmode C\!\!\!\!I\else$C\!\!\!\!I$\fi}

\def\K{\relax\ifmmode I\!\!K\else$I\!\!K$\fi}

\def\N{\relax\ifmmode I\!\!N\else$I\!\!N$\fi}

\newcounter{defcounter}[section]
\newenvironment{definition}%
{\vspace{0.1cm}\begin{sloppypar}\noindent\stepcounter{defcounter}{\bfseries
Definition
      \thesection.\thedefcounter}}%
{\end{sloppypar}\vspace{0.1cm}}

\newtheorem{lemma}{Lemma}[section]
\newtheorem{theorem}{Theorem}[section]

\newtheorem{proposition}{Proposition}[section]

\newcommand{\proof}{{\bf Proof.} }

\newcommand{\qed}{\hfill $\square$}

\begin{document}
\thispagestyle{empty}
\begin{center}
{\Large {\bf Topological chaos in the two-gene Andrecut-Kauffman model}}
\end{center}
\begin{center}J\"org Neunh\"auserer\\
Technical University of Braunschweig \\
joerg.neunhaeuserer@web.de
\end{center}
\begin{center}
\begin{abstract}
The two-gene Andrecut-Kauffman model describes gene expression within cells by a two-dimensional discrete dynamical system.
Using periodic orbits of an one dimensional subsystem we, rigourously prove the existence of topological chaos in this system. In addition we compute parameter domains of the system where we have topological chaos.\\
{\bf MSC 2020: Primary: 37N25, 65P20 Secondary: 37B40, 92C42}~\\
{\bf Key-words: gene expression, two-dimensional discrete dynamical system, one-dimensional dynamics, topological chaos}
\end{abstract}
\end{center}
\section{Introduction}
To understand the dynamics of biological processes within cells at the molecular level, gene expression models are essential.
Mathematical models that accurately reflect these processes allows to predictic the behavior of the biological system. This even contributes
to the development and advancement of therapeutic strategies diseases. 
Gene expression models describe the complex networks of interactions between genes,
transcripts, and proteins, and leads to a better understanding of cellular responses to stress and cell differentiation processes.
The expression of a single gene does not describe the complex
processes occurring within a living cell.  When we consider systems with more genes, their corresponding proteins can act as transcription factors,
binding to different promoters, which leads to a network of dependencies. 
In the literature we find many attempts to accurately model gene expression networks.
Classical models like the Michaelis-Menten
model \cite{[MM]} or the Goodwin’s model \cite{[GO]}
use systems of differential equations. This models describe the concentrations
of gene products, such as mRNA and proteins, in a smooth manner over time.
These systems form the basis for analyzing complex dynamic dependencies among components of
a biological system and have been further developed by many authors. The limitation of such models is the assumption of an immediate response of the
system to changes in parameters.  More realistic are time delay models, which pay attention to the time delay between the initiation of transcription and
the protein synthesis. In the resent works \cite{[BR1],[BR2],[BR3]} the authors
analyze time continuous mathematical models of Hes1 gene expression, incorporating various aspects such
as time delays, feedback mechanisms and oscillations. Similarly, the
authors of \cite{[BR4]} investigate a generalized p53-Mdm2 model and demonstrate that
oscillations in gene expression may arise.\\
The Andrecut-Kauffman model we study in this paper has its origin in the work of Kauffman \cite{[KA]} who, used Boolean networks to describe gene expression. 
Andrecut \cite{[AN]} considered random Boolean networks as a model and used mean field theory to calculate
probabilities of node states in such networks. Then Andrecut and Kauffman \cite{[AK2]} proposed a model in which individual
reactions between genes have been replaced with their average results, making the outcome deterministic.
The discrete-time model we study here was proposed one year later in \cite{[AK]}. The
number of genes was fixed at 2 and only homo-multimers of length exactly $n$ were
taken into consideration. This led to a two dimensional discrete dynamical system, which seems to represent the
complex dynamics of the original model as well. We are interested in this model from a mathematical perspective since Andrecut and Kauffman provide numerical evidence for a chaotic dynamics of the system.\\
The rest of the paper is organized as follows. In the next section we introduce the Andrecut-Kauffman model and given an interpretation of the variables and parameters of the system. Then in section 4 we introduce the notation from topological dynamics we use. We cite two results on topological chaos that form the basis of our argumentation in the next section. Section 5 contains the main result of the paper and its proof. In section 6 we apply this result. For a specific choice of parameter values we prove that the system is topological chaotic. Moreover we compute parameter domains for which we have topological chaos. The paper concludes with a discussion of our results and suggestions for future research.
\section{The Andrecut-Kauffman model}
For parameters $m\in\mathbb{N}$, $\epsilon\in[0,1]$,  $a>0$ and $b\in(0,1)$ we consider the two-dimensional time-discrete Andrecut-Kauffman model, given by
\[ T_{m,\epsilon,a,b}:\mathbb{R}^{2}_{\ge 0} \to \mathbb{R}^{2}_{\ge 0}\] with
\[ T_{m,\epsilon,a,b}(x,y)=\left(bx+\frac{a}{1+(1-\epsilon)x^m+\epsilon y^m},by+\frac{a}{1+\epsilon x^m+(1-\epsilon)y^m}\right).  \]
Note that for $(x,y)\in [0,\frac{a}{1-b}]^2$ we have 
\[ 0<bx+\frac{a}{1+(1-\epsilon)x^m+\epsilon y^m}<bx+a\le\frac{a}{1-b}\]
and
\[ 0<by+\frac{a}{1+\epsilon x^m+(1-\epsilon)y^m}<by+a\le \frac{a}{1-b},\]
hence
$T_{m,\epsilon,a,b}(x,y)\in [0,\frac{a}{1-b}]^2$. We may restrict our attention to this compact square $[0,\frac{a}{1-b}]^2$. \\
According to \cite{[AK]} the variables  $x$ and $x$ represent the concentration levels of transcription factor proteins
for two genes x and y. The orbit
\[ \mathcal{O}(x,y)=\{T_{m,\epsilon,a,b}^{n}(x,y)|n\ge 0\}\]
describes the evolution of these concentrations for the initial value $(x,y)$.
The parameters $a,b$ represent properties of the chemical reactions, describing the 
complex steps of gene expression. The parameter $a$ reflects
the rate of combined transcription and translation over a fixed time interval. The
parameter $b$ depends on the rate of protein degradation for a fixed time
interval. $\epsilon$ is a parameter that describes the coupling between the genes, where it is assumed that the coupling is symmetric. The exponent $n$ denotes the number
of monomers of a given protein, that are subject to the multimerization reaction.\\
In the case $m = 1$ and $m=2$ the model seems to have non chaotic behavior, there are only stable fixpoints and periodic orbits. We are thus interested in the case  $m\ge 3$.  In the case $m = 3$ and $m=4$ Andrecut and Kauffman \cite{[AK]} approximated lokal Lyapunov-exponents of the system for specific orbits, given various parameters and initial values, which are not explicitly stated. If the calculation indicates that a lokal Lyapunov-exponent is positive, they conclude that the dynamics of the system, with the parameters considered, is chaotic. From a mathematical point of view this kind of reasoning is not satisfactory. In fact a positive Lyapunov-exponent with respect to an absolutely continuous ergodic measure or an ergodic measure on a fractal invariant set implies topological chaos. Such a measure has positive metric entropy which implies positive topological entropy of the system, see \cite{[KH]}. But it is not clear and difficult to prove that such a measure exists for the Andrecut-Kauffman model. Even if we assume the existence of such a measure, it is possible that we choose non-typical initial values to approximate a local Lyapunov-exponent. The set of such 
initial values may be large, see \cite{[BS]}. In this case the calculations of Lyapunov-exponents does not show that the dynamics of the system is chaotic.\\  
We develop here another approach to proof topological chaos in Andrecut-Kauffman model 
and to find parameter domains where the dynamics of the system is chaotic.  We show the existence of periodic orbits with an odd period of 
an one dimensional subsystem of the Andrecut-Kauffman system. This guarantees positive topological entropy of the subsystem and of the system itself, which is topological chaos. For a specific chose parameter values we thus get a rigourous of proof of chaos in Andrecut-Kauffman model. 
\section{Topological entropy and chaos}
In this section we recall the definition of topological entropy and topological chaos. We will cite two results, which will be used in the sequel.\\
Let $X$ be a metric space and ${\mathfrak U}$ be an open covering of $X$.
The entropy of this covering is $H({\mathfrak U})=\log(\sharp{\mathfrak U})$, where $\sharp{\mathfrak U}$ is the minimal number of sets in ${\mathfrak U}$ needed to cover $X$. For two coverings ${\mathfrak U}_{1},{\mathfrak U}_{2}$ auf $X$ we define the common refinement by  
\[  {\mathfrak U}_{1}\vee {\mathfrak U}_{2}=\{O_{1}\cap O_{2}|O_{1}\in {\mathfrak U_{1}}, O_{2}\in {\mathfrak U_{2}}\}.\]
A refinement of a covering obviously increases the entropy of the covering. We now have the necessary notations to define topological entropy and topological chaos.   
\begin{definition} Let $X$ be a compact metric space and let $T:X\to X$ be continuous. The entropy of the system $(X,T)$ with respect to a covering ${\mathfrak U}$ of $X$ is
\[ h(T,{\mathfrak U})=\lim_{n\to\infty}\frac{1}{n}H({\mathfrak U}\vee T^{-1}({\mathfrak U})\vee\dots\vee T^{-n+1}({\mathfrak U}))
\footnote{The limit exsist, since $a_{n}=H({\mathfrak U}\vee T^{-1}({\mathfrak U})\vee\dots\vee T^{-n+1}({\mathfrak U}))$ is a subadditive sequence, 
that means $a_{n+m}\le a_{n}+a_{m}$.}\]
and the entropy of the system is
\[ h(T)=\sup\{h(T,{\mathfrak U})~|~{\mathfrak U}\mbox{ ist an open covering of }X\}.\]
$(X,T)$ topological chaotic, if $h(T)>0$.
\end{definition} 
~\\The definition we present here goes back to Adler, Konheim and McAndrew  \cite{[AD]}. Another approach, which is a more involved but allows the definition of entropy for non compact systems can be found in Bowen \cite{[BOW]}. The topological entropy measure the complexity and unpredictability of a topological dynamical system and thus provide a quantitative definition of chaos.\\
The following proposition is well known and follows directly from the definition: 
\begin{proposition}
Let $X$ be a compact metric space and let $T:X\to X$ be continuous. If $K\subset X$ is a compact $T$-invariant set we have
\[ h(T_{|K})\le h(T).\]
Especially if the subsystem $(K,T)$ is topological chaotic the system $(X,T)$ is topological chaotic as well. 
\end{proposition}
For continuous interval maps the theorem of Misiurewicz \cite{[MI]} relates topological chaos to the existence of certain periodic orbits:  
\begin{theorem}
Let $I$ be a closed interval and let $T:I\to I$ be continuous. $(I,T)$ is topological chaotic if and only if there exists a periodic orbit with a period that is not a power of $2$.
\end{theorem}
Readers interested in interval dynamics should consider the nice book of Ruette \cite{[RUE]}, where one finds Misiurewicz theorem as theorem 4.48. 
\section{The main result}
For parameters $m\in\mathbb{N}$, $a>0$ and $b\in(0,1)$ we consider the one dimensional map $T_{m,a,b}(x):[0,\frac{a}{1-b}]\to [0,\frac{a}{1-b}]$, given by
\[T_{m,a,b}(x)=bx+\frac{a}{1+x^m}.\] 
Note that $T_{m,a,b}$ describes the action of the Andrecut-Kauffman map $T_{m,\epsilon,a,b}$ from section 2 on the diagonal \[ D=\{(x,x)|x\in[0,\frac{a}{1-b}]\}\subset [0,\frac{a}{1-b}]^2\]\\ 
for all $\epsilon \in[0,1]$.\\
The following lemma is useful when we want to find periodic orbits of the map $T_{m,a,b}$, that are not fixpoints.
\begin{lemma}
The map $T_{m,a,b}(x):[0,\frac{a}{1-b}]\to [0,\frac{a}{1-b}]$ has unique fixpint $x_{m,a,b}\in (0,\frac{a}{1-b})$. 
\end{lemma}
\proof Let 
\[ h(x)=T_{m,a,b}(x)-x=-(1-b)x+\frac{a}{1+x^m}. \]
We have $h(0)=a>0$ and $h(a/(1-b))<0$. Moreover $h$ is strictly monoton decreasing since
\[ h^{\prime}(x)=-(1-b)-\frac{amx^{m-1}}{(1+x^{m})^2}<0\]
for$x\ge 0$. It follows that there is a unique $x_{m,a,b}\in (0,\frac{a}{1-b})$ with $h(x_{m,a,b})=0$, which is equivalent to $T_{m,a,b}(x_{m,a,b})=x_{m,a,b}$.  
\qed ~\\~\\
Now we are prepared to proof the main theorem of this paper:
\begin{theorem}
Let $m\in\mathbb{N}$, $a>0$ and $b\in(0,1)$. If there exists a $x\in (0,\frac{a}{1-b})\backslash \{x_{m,a,b}\}$ and an odd numbers $n\in\mathbb{N}$ with $T_{m,a,b}^{n}(x)=x$, than the Andrecut-Kauffman model $([0,\frac{a}{1-b}]^2,T_{m,\epsilon,a,b})$ is topological chaotic for all $\epsilon\in[0,1]$.
\end{theorem}
\proof Assume that $T_{m,a,b}^{n}(x)=x$ with $x\in (0,\frac{a}{1-b})\backslash \{x_{m,a,b}\}$ for an odd number $n$. Since $x\not=x_{m,a,b}$
$x$ is by lemma 4.1 not the fixpoint of the system.
Since $n$ is odd it follows from theorem 3.1 that the system $([0,\frac{a}{1-b}],T_{m,a,b})$ is topological chaotic. Hence the system $(D,T_{m,a,b}\times T_{m,a,b})$ is topological chaotic. But this is a subsystem of $([0,\frac{a}{1-b}]^2,T_{m,\epsilon,a,b})$ for all $\epsilon\in[0,1]$. By proposition 3.1 this system is thus topological chaotic as well.    
\qed\\\\
If $T_{m,a,b}^{n}(x)=x$ holds for $x\in (0,\frac{a}{1-b})\backslash \{x_{n,\epsilon,a,b}\}$ and an even numbers $n$ one would have to check in addition that $x$ has not a period that is power of $2$, to guarantee topological chaos of the system. We will not apply this variant of theorem 4.1 in the next section.\\
\section{Applications}
We first consider three concrete applications of theorem 4.1. For $m=3$, $a=50$ and $b=0.1$ the map $T_{3,50,0.1}:[0,500/9]\to [0,500/9]$ is 
given by   
\[T_{3,50,0.1}(x)=\left(0.1x+\frac{50}{1+x^3}\right).\]
We do not find a periodic Orbit of period $3$ for this map but we find a periodic Orbit of period $5$.
The Orbit of $2$ of the map may easily be calculated even without using a machine
\[ \mathcal{O}(2)=\{T_{3,50,0.1}^{n}(2)| 0\le n\le 5\}\]                                  
\[ \approx\{2, 5.75556, 0.836433, 31.6257, 3.16415, 1.84645\}.\]
We have $T_{3,50,0.1}^{5}(2)<2$ hence by continuity there is an $p\in (0,2)$ with $T_{3,50,0.1}^{5}(p)=p$. Since $T_{3,50,0.1}(2)>2$ we have $T_{3,50,0.1}^{5}(p)>p$ hence $p$ is not the fixpoint of the map. By theorem 4.1 we conclude that the system $([0,\frac{500}{90}]^2,T_{3,\epsilon,50,0.1})$ is topological chaotic for all $\epsilon\in[0,1]$. \\
Now we keep $a=50$ and $b=0.1$ but choose the exponent $m=4$. In this case we have
\[ \mathcal{O}(2)=\{T_{4,50,0.1}^{n}(2)| 0\le n\le 5\}\]                                  
\[ \approx\{2,3.14118, 0.822467, 34.3855, 3.43858, 0.698963 \}.\]
As above we get from theorem 4.1 that $([0,\frac{500}{90}]^2,T_{4,\epsilon,50,0.1})$ is topological chaotic for all $\epsilon\in[0,1]$.\\
For $m=5$, $a=50$ and $b=0.1$ we consider the Orbit of $3/2$. Since
 \[ \mathcal{O}(3/2)=\{T_{5,50,0.1}^{n}(3/2)| 0\le n\le 5\}\]  
\[\approx\{3/2, 5.96818, 0.603421, 46.3566, 4.63566, 0.486912\}\] 
the system $([0,\frac{500}{90}]^2,T_{4,\epsilon,50,0.1})$ is topological chaotic for all $\epsilon\in[0,1]$. By continuity these results remain true for parameters $(a,b)$ in some neighborhood of $(50,0.1)$. We claim that these applications of theorem 4.1 provide the first rigourous proof of chaos in the Andrecut-Kauffman model.\\
Now we consider the set
\[ \mathcal{C}_{n}(m):=\{(a,b)\in (0,100)\times (0,0.5)| \exists x\in (0,\frac{a}{1-b})\backslash \{x_{m,a,b}\}~:~T_{m,a,b}^{n}(x)=x\}\]
for $m\ge 3$ and $n$ odd. By theorem 4.1 these are parameter domains for which the Andrecut-Kauffman system 
$([0,\frac{a}{1-b}]^2,T_{m,\epsilon,a,b})$ ist topological chaotic for all $\epsilon\in[0,1]$.  \\
If $n$ and $m$ are small enough we can use a computer algebra system to represent $T_{m,a,b}^{n}(x)$ symbolically and an equation solver und find
appropriate solutions of $T_{m,a,b}^{n}(x)=x$ for given parameters $a,b$ and $m$. In this way we approximate the set $\mathcal{C}_{n}(m)$ for $500\times500$ parameters in $(0,100)\times (0,0.5)$, where we consider the exponents $m=3,4,5$ of the model. We find no periodic orbits of period $3$ thus $\mathcal{C}_{3}(m)$ seems to be empty. In figure 1
we display an approximation of $\mathcal{C}_{5}(3)$ and $\mathcal{C}_{7}(3)$ representing the domains where we find periodic orbits of period $5$ and $7$ for the exponent $m=3$. 
\begin{figure}
\vspace{0pt}\hspace{0pt}\scalebox{0.37}{\includegraphics
{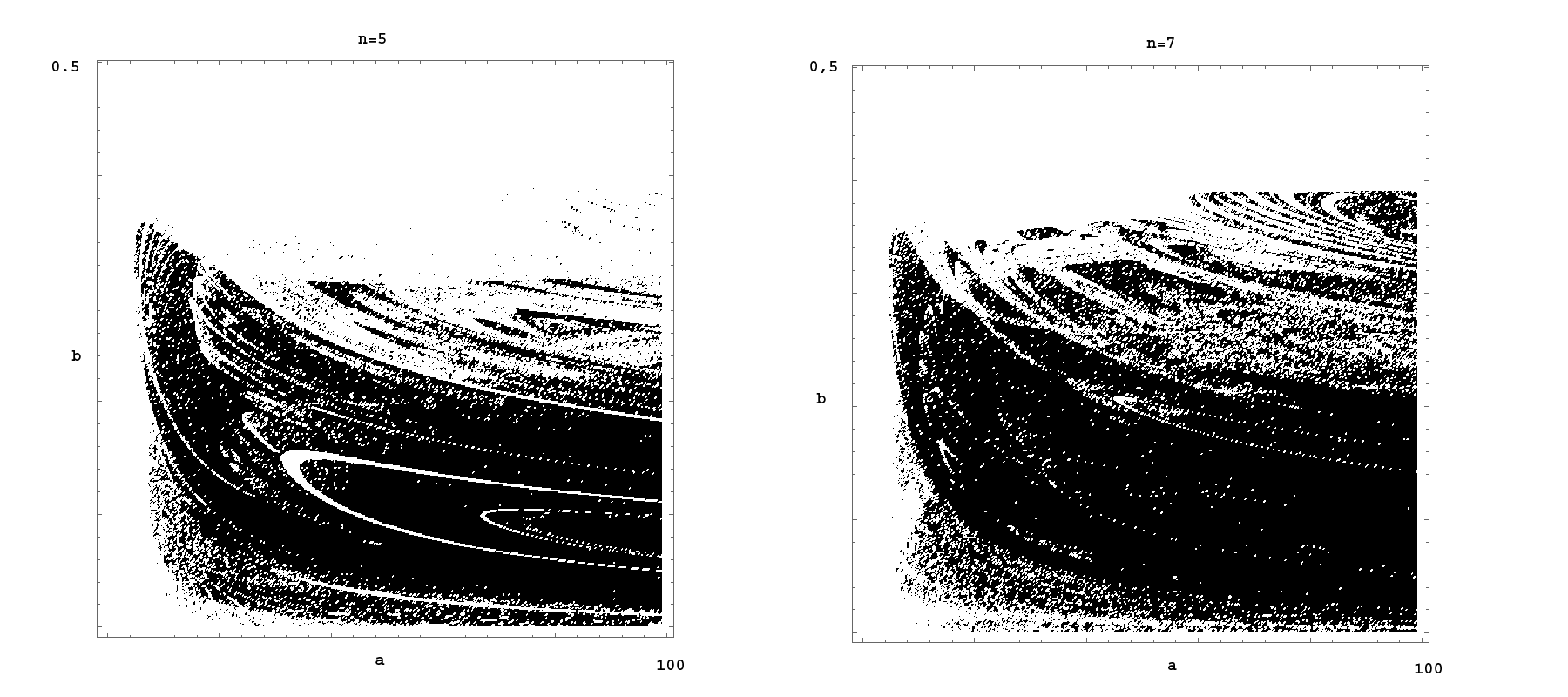}}
\caption{The domains $\mathcal{C}_{5}(3)$ and $\mathcal{C}_{7}(3)$}
\end{figure}
In figure 2 we display an approximation of $\mathcal{C}_{5}(4)$ and $\mathcal{C}_{7}(4)$ and in figure 3 we display an approximation of $\mathcal{C}_{5}(5)$ and $\mathcal{C}_{7}(5)$, which are domains of chaos for Andrecut-Kauffman model in the case of exponent $m=4$ resp. $m=5$. Of course these figures only give numerical evidence. To get rigourous results one can check specific pairs of parameters in the way we did above. 
\begin{figure}
\vspace{0pt}\hspace{0pt}\scalebox{0.37}{\includegraphics
{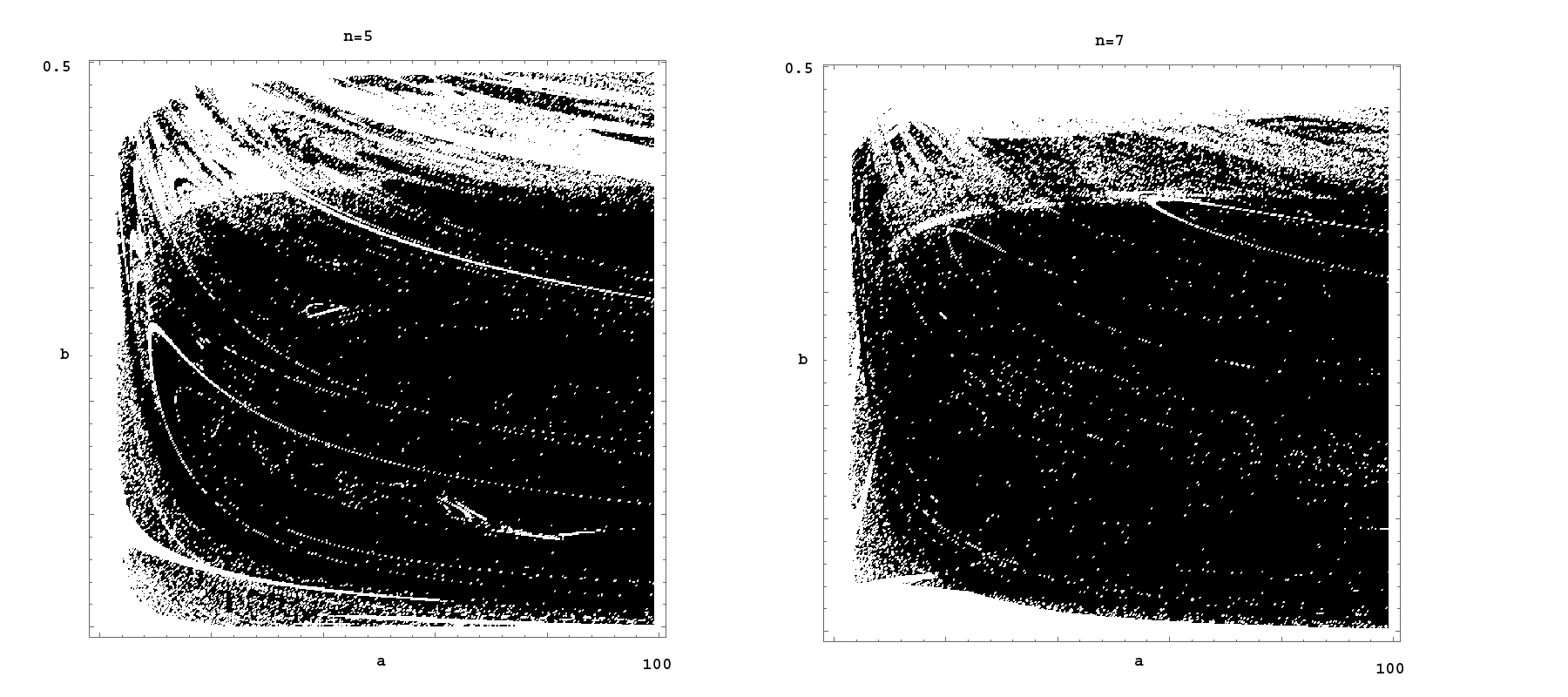}}
\caption{The domains $\mathcal{C}_{5}(4)$ and $\mathcal{C}_{7}(4)$}
\end{figure}
\begin{figure}
\vspace{0pt}\hspace{0pt}\scalebox{0.37}{\includegraphics
{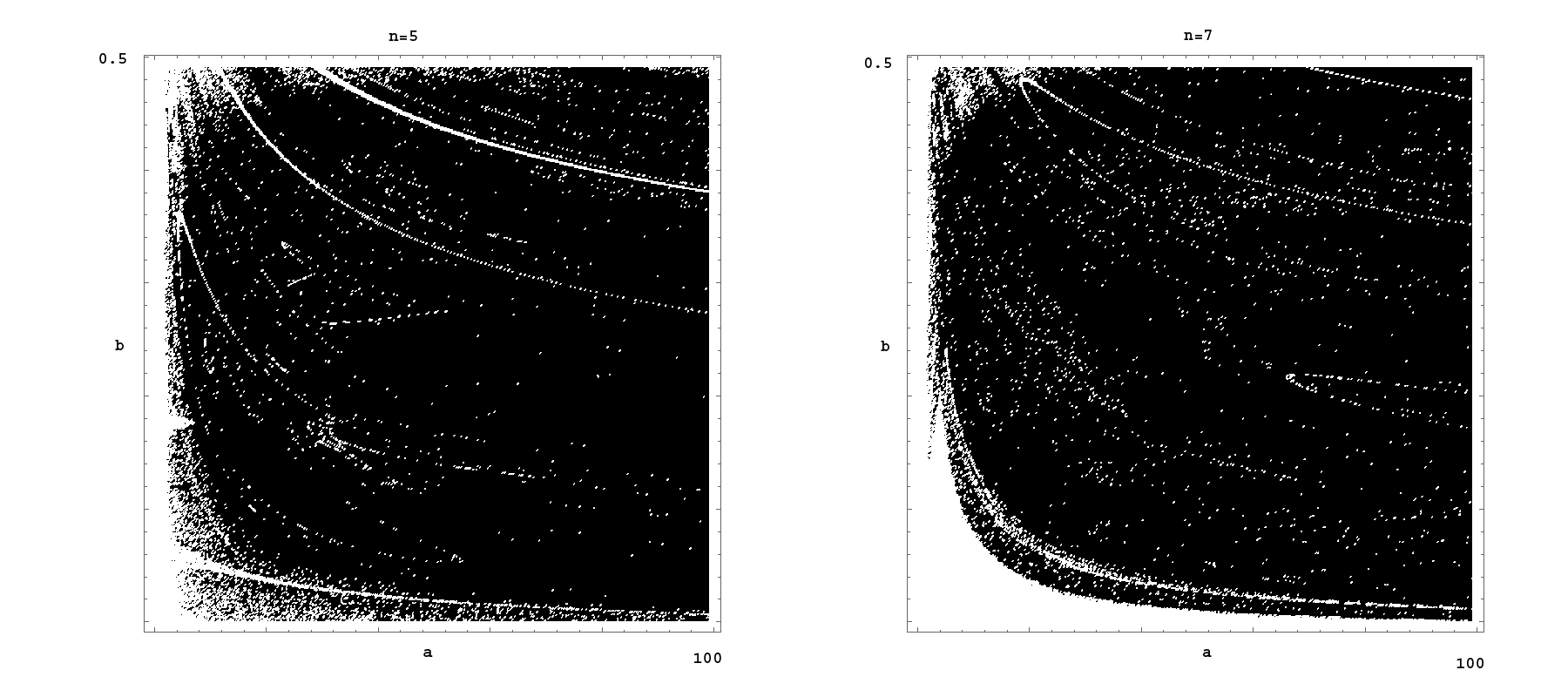}}
\caption{The domains  $\mathcal{C}_{5}(5)$ and $\mathcal{C}_{7}(5)$}
\end{figure}
\section{Conclusion and future research directions}
We have seen that using periodic orbits of an one dimensional subsystem it is possible to rigorously prove the existence of topological chaos in the two-gene Andrecut-Kauffman model for gene expression. With this approach it is also possible to compute parameter domains of the system, where we have topological chaos. Our results follow from the meanwhile classical theory of topological dynamical systems. It would be interesting to introduce three-gene or more generally $n$-gen models of gene expression similar to the model of Andrecut and Kaufmann. Such models would allow a more realistic description of cell dynamics at the molecular level. We suppose that we could find chaotic one-dimension subsystems in such models, which allow a proof that the systems are topological chaotic for a specific choice parameters. 

\end{document}